\documentclass{article}

\usepackage{proof}
\usepackage{latexsym}
\usepackage{amsfonts}

\newcommand{\bremark}{\begin{remark}}
\newcommand{\eremark}{\end{remark}}

\newcommand{\blem}{\begin{lemma}}
\newcommand{\elem}{\end{lemma}}
\newcommand{\bth}{\begin{theorem}}
\newcommand{\eth}{\end{theorem}}
\newcommand{\benu}{\begin{enumerate}}
\newcommand{\eenu}{\end{enumerate}}
\newcommand{\bdes}{\begin{description}}
\newcommand{\edes}{\end{description}}
\newcommand{\bdf}{\begin{definition}}
\newcommand{\edf}{\end{definition}}
\newcommand{\bcor}{\begin{cor}}
\newcommand{\ecor}{\end{cor}}
\newcommand{\bprp}{\begin{proposition}}
\newcommand{\eprp}{\end{proposition}}
\newcommand{\bmlem}{\begin{mlemma}}
\newcommand{\emlem}{\end{mlemma}}
\newcommand{\bclm}{\begin{claim}}
\newcommand{\eclm}{\end{claim}}
\newcommand{\bprf}{{\bf Proof}.\hspace{2mm}}

\newcommand{\eprf}{\hspace*{\fill} $\Box$}

\newcommand{\eprfclm}{\hspace*{\fill} {\it End of Proof of Claim}}

\newcommand{\beqn}{\begin{equation}}
\newcommand{\eeqn}{\end{equation}}
\newcommand{\beqnarr}{\begin{eqnarray}}
\newcommand{\eeqnarr}{\end{eqnarray}}
\newcommand{\beqnarrs}{\begin{eqnarray*}}
\newcommand{\eeqnarrs}{\end{eqnarray*}}
\newcommand{\smlskp}{\\ \smallskip \\}
\newcommand{\spand}{\,\&\,}

\newtheorem{theorem}{Theorem}[section]
\newtheorem{definition}[theorem]{Definition}
\newtheorem{proposition}[theorem]{Proposition}
\newtheorem{lemma}[theorem]{Lemma}
\newtheorem{cor}[theorem]{Corollary}
\newtheorem{mlemma}{Main Lemma}
\newtheorem{claim}{Claim}

\newtheorem{remark}[theorem]{Remark}

\newcommand{\alp}{\alpha}

\newcommand{\veps}{\varepsilon}
\newcommand{\del}{\delta}
\newcommand{\Del}{\Delta}
\newcommand{\ome}{\omega}
\newcommand{\Ome}{\Omega}
\newcommand{\bet}{\beta}
\newcommand{\gam}{\gamma}
\newcommand{\Gam}{\Gamma}
\newcommand{\kap}{\kappa}
\newcommand{\sig}{\sigma}
\newcommand{\Sig}{\Sigma}

\newcommand{\lam}{\lambda}
\newcommand{\Lam}{\Lambda}
\newcommand{\vphi}{\varphi}

\newcommand{\fal}{\forall}
\newcommand{\exi}{\exists}
\newcommand{\rarw }{\rightarrow}
\newcommand{\Rarw }{\Rightarrow}

\newcommand{\lrarw}{\leftrightarrow}
\newcommand{\Lrarw}{\Leftrightarrow}

\newcommand{\calL}{{\cal L}}

\newcommand{\incl}{\subseteq}

\title{A sneak preview of proof theory of ordinals\thanks{This is a revised version of the r\'esum\'e for a talk at Kobe seminar on Logic and Computer Science
, 5-6 Dec.1997}}
\author{Toshiyasu Arai
 \thanks{I would like to thank Prof. Y. Kakuda and Dr. M. Kikuchi for hospitality during my visit to Kobe.}
 \\
Faculty of Integrated Arts and Sciences\\
Hiroshima University\\
Higashi-Hiroshima, 739 Japan\\
arai@mis.hiroshima-u.ac.jp
\thanks{
current address:
Graduate School of Science,
Chiba University,
1-33, Yayoi-cho, Inage-ku,
Chiba, 263-8522, JAPAN,
tosarai@faculty.chiba-u.jp
}
}
\date{Dec. 5 1997}
\begin{document}
\maketitle
\begin{abstract}
This talk is a sneak preview of the project,
'proof theory for theories of ordinals'. Background,
aims, survey and furture works on the project are given.
Subsystems of second order arithmetic are embedded in
recursively large ordinals and then the latter are analysed.
We scarcely touch upon proof theoretical matters.
\end{abstract}

\section {Proof theory \`a la Gentzen-Takeuti}\label{sec:1}
Let T be a sound and recursive theory containing arithmetic. The {\it proof-theoretical ordinal\/} 
$|\mbox{T}|_{\Pi^{1}_{1}}<\ome^{CK}_{1}$ is defined by the ordinal:
\[\sup\{\alp<\ome^{CK}_{1}: \mbox{T}\vdash Wo[\prec_{\alp}] \mbox{ for some recursive well ordering } \prec_{\alp} \mbox { of type } \alp\}\]
($Wo[\prec]$ denotes a $\Pi^{1}_{1}$-sentence saying that $\prec$ is a well ordering.)
\benu
\item (Gentzen 1936, 1938, 1943) $|\mbox{PA}|_{\Pi^{1}_{1}}=|\mbox{ACA}_{0}|_{\Pi^{1}_{1}}=\veps_{0}$
\item (Takeuti 1967) $|\Pi^{1}_{1}\mbox{-CA}_{0}|_{\Pi^{1}_{1}}=|O(\ome,1)|_{<_{0}}=\psi_{\Ome}\Ome_{\ome}$ and\\
$|\Pi^{1}_{1}\mbox{-CA+BI}|_{\Pi^{1}_{1}}=|O(\ome+1,1)|_{<_{0}}=\psi_{\Ome}\veps_{\Ome_{\ome}+1}$
\eenu
{\bf Axiom schemata in second order arithmetic}. Let $\Phi$ denote a set of formulae in the language of second order arithmetic.
\benu
\item $\Phi\mbox{-CA}$: For each $\vphi\in\Phi$
\[\fal Y\exi X[X=\{n\in\ome: \vphi(n,Y)\}]\]
\item $\Phi^{-}\mbox{-CA}$ denotes the set-parameter free version of $\Phi\mbox{-CA}$:
\[\exi X[X=\{n\in\ome: \vphi(n)\}]\]
\item $\Del^{1}_{n}\mbox{-CA}$: For $\vphi,\psi\in\Sig^{1}_{n}$
\[\{n\in\ome:\vphi(n,Y)\}=\{n\in\ome:\neg\psi(n,Y)\}\rarw \exi X[X=\{n\in\ome: \vphi(n,Y)\}]\]
\item $\Phi\mbox{-AC}$: For each $\vphi\in\Phi$
\[\fal n\exi X \vphi(n,X)\rarw \exi \{X_{n}\}\fal n\vphi(n,X_{n})\]
\item $\Phi\mbox{-DC}$: For each $\vphi\in\Phi$
\[\fal n\fal X\exi Y\vphi(n,X,Y)\rarw \exi\{X_{n}\}\fal n\vphi(n,X_{n},X_{n+1})\]
\item BI: For each formula $\vphi$
\[Wf[X]\rarw TI[X,\vphi]\]
\eenu

 Proof theory \`a la Gentzen-Takeuti \cite{G3}, \cite{T} proceeds as follows;
\bdes
\item[(G1)] Let $P$ be a proof whose endsequent $\Gam$ has a restricted form, e.g., an arithmetical sequent. Define a reduction procedure $r$ which rewrites such a proof $P$ to yield another proofs $\{r(P,n):n\in I\}$ of sequents $\Gam_{n}$  provided that $P$ has not yet reduced to a certain canonical form. 

For example when we want to show that the arithmetical sequent $\Gam$ is true, the sequents $\Gam_{n}$ are chosen so that $\Gam$ is true iff every $\Gam_{n}\, (n\in I)$ is true. Also if $P$ is in an irreducible form, then the endsequent is true outright.  
\item[(G2)] From the structure of the proof $P$, we abstract a structure related to this procedure $r$ and throw irrelevant residue away. Thus we get a finite figure $o(P)$.

We call the figure $o(P)$ the {\em ordinal diagram} (o.d.) following G. Takeuti \cite{T}. Let ${\cal O}$ denote the set of o.d.'s.
\item[(G3)] Define a relation $<$ on ${\cal O}$ so that $o(r(P,n))<o(P)$ for any $n\in I$.
\item[(G4)] Show the relation $<$ on ${\cal O}$ to be well founded. \\
Usually $<$ is a linear ordering and hence $({\cal O},<)$ is a notation system for ordinals.
\edes
When the endsequent of a proof $P$ is an arithmetical sequent, we in fact construct an $\ome$ cut-free proof of the sequent whose height is less than or equal to (the order type of) the o.d. $o(P)$ attached to $P$.

O.d.'s are constructed so that each constuctor for o.d.'s reflects a reduction step on proofs.

We attach an o.d. $o(\Gam;P)$ to each sequent $\Gam$ occurring in a proof $P$. The o.d. $o(\Gam;P)$ is built by applying constructors for o.d.'s. Applied constructors in building the term $o(\Gam;P)$ correspond to the inference rules occurring above $\Gam$.

\section{$\ome$-proofs}\label{sec:2}
In the latter half of 60's Sch\"utte, Tait, Feferman et.al analysed predicative parts of second order arithemetic using infinitary proofs with $\ome$-rule: infer $\fal n A(n)$ from $A(n)$ for any $n\in\ome$. Their main result is
\[|\mbox{ATR}_{0}|_{\Pi^{1}_{1}}=\Gam_{0}=_{df}\min\{\alp>0:\fal\bet,\gam<\alp(\vphi\bet\gam<\alp)\}\]
where $\vphi$ denotes the binary Veblen function: For each $\alp<\ome_{1}$ define inductively a normal (strictly increasing and continuous) function $\vphi_{\alp}:\ome_{1}\rarw\ome_{1}$ as follows: First set $\vphi 0\bet=\vphi_{0}\bet=\ome^{\bet}$. Since the ranges $rng(\vphi_{\bet})$ of $\vphi_{\bet}\, (\bet<\alp)$ are club sets in $\ome_{1}$, so are their fixed points $fp(\vphi_{\bet})=\{\gam<\ome_{1}:\vphi_{\bet}\gam=\gam\}$. Thus the intersection $\bigcap\{fp(\vphi_{\bet}):\bet<\alp\}$ is also a club set in $\ome_{1}$. $\vphi_{\alp}$ is defined to be the enumerating function of the set $\bigcap\{fp(\vphi_{\bet}):\bet<\alp\}$.

\section{Buchholz-Pohlers}\label{sec:3}
In their Habilitationsscriften (1977) Buchholz \cite{Buchholz77}($\Ome_{\mu+1}$-rule) and \\
Pohlers \cite{Pohlers77}(local predicativity method) analysed theories for iterated inductive definitions. These theories formalize least fixed points of positive elemenary induction on $\ome$. For a monotone operator $\Gam:{\cal P}(\ome)\rarw{\cal P}(\ome)$ define inductively sets $I^{\Gam}_{\alp}$ by
\[I^{\Gam}_{\alp}=\Gam(I^{\Gam}_{<\alp})\mbox{ with } I^{\Gam}_{<\alp}=\bigcup\{I^{\Gam}_{\bet}:\bet<\alp\}\]
$I^{\Gam}=\bigcup\{I^{\Gam}_{\alp}:\alp<\ome^{CK}_{1}\}$ is the least fixed point of $\Gam$.

When $\Gam$ is given by a positive elementary formula $A(X^{+},n)$, $\Gam(X)=\{n\in\ome:\ome\models A[X,n]\}$, we write $I^{A}$ for $I^{\Gam}$ and $I^{A}_{\alp}$ for $I^{\Gam}_{\alp}$.

Iterated $\Pi^{1-}_{1}\mbox{-CA}_{\nu}$ can be simulated in $\mbox{ID}_{\nu}$ since
\[\Pi^{1}_{1} \mbox{ on } \ome=\mbox{ inductive on } \ome=\Sig_{1} \mbox{ on } \ome^{CK}_{1}\]
This is seen from Brower-Kleene $\Pi^{1}_{1}$-normal form: for each $\Pi^{1}_{1}$ $A(n,X)\, (n\in\ome,X\incl \ome)$ there exists a recursive relation $Q_{n,X}$ on $\ome$, i.e., there exists a recursive function $\{e\}^{X}(n)$ so that
\beqnarr
&& A(n,X)\lrarw Wf[Q_{n,X}]\lrarw \{e\}^{X}(n)\in{\cal O}\lrarw \nonumber \\
&& \exi f\in L_{\ome^{CK}_{1}}[f \mbox{ collapsing function for } Q_{n,X}]\lrarw \nonumber \\
&& \exi\alp<\ome^{CK}_{1}\fal x\in dom(Q_{n,X})(x\in I^{Q_{n,X}}_{\alp}) \label{eqn:BK}
\eeqnarr
($x\in I^{Q_{n,X}}_{\alp}$ designates the order type of $Q_{n,X}| x$ is less than or equal to $\alp$.)

They showed ({\it cf}. \cite{LNM897} (1981).)
\benu
\item $|\mbox{ID}_{\nu}|_{\Pi^{1}_{1}}=|\Pi^{1-}_{1}\mbox{-CA}_{\nu}|_{\Pi^{1}_{1}}=\psi_{\Ome}\veps_{\Ome_{\nu}+1}$
\item $|\mbox{ID}_{<\lam}|_{\Pi^{1}_{1}}=|\Pi^{1-}_{1}\mbox{-CA}_{<\lam}|_{\Pi^{1}_{1}}=\psi_{\Ome}\Ome_{\lam}$ for limit $\lam$
\eenu
$\Ome_{\nu}$ denotes either $\ome_{\nu}$ or the continuous closure of the enumerating function of the recursively regular ordinals.

\bremark. 
{\rm Recently (May 1997) Buchholz \cite{Buchholz 1997} shows that Sch\"utte's cut eliminaion procedure for infinitary proofs with $\ome$-rule is nothing but the infinitary image of Gentzen's, $\mbox{Gentzen}^{\infty}=\mbox{Sch\"utte}$ and $\mbox{Takeuti}^{\infty}=\mbox{Buchholz}$.

I conjecture that $\mbox{Arai}^{\infty}=\mbox{Pohlers-J\"ager}$ for $\mbox{KP}\ome$.}
\eremark

\section{J\"ager}\label{sec:4}
   G. J\"ager \cite{J} has shifted an object of proof-theoretic study to set theories from second order arithmetic.

\begin{definition}($\Pi_2^\Ome$-ordinal of a theory)  
{\rm Let T be a recursive theory of sets such that}
 $\mbox{{\rm KP}}\omega\subseteq \mbox{{\rm T}}\subseteq \mbox{{\rm ZF+V=L}}$, {\rm where} $\mbox{{\rm KP}}\omega$ {\rm denotes Kripke-Platek set theory with the Axiom of Infinity.} {\rm For a sentence} $A$ {\rm let}  $A^{L_{\alpha}}$ {\rm denote the result of replacing unbounded quantifiers} $Qx \, (Q\in\{ \forall, \exists\})$ {\rm in} $A$ {\rm by} $Qx\in L_{\alpha}$. {\rm Here for an ordinal} $\alpha\in Ord$ $L_\alpha$  {\rm denotes an intial segment of G\"odel's constructible sets.} {\rm Let} $\Omega$ {\rm denote the (individual constant corresponding to the) ordinal} $\omega^{CK}_1$. {\rm If} $T\not\vdash\exists \omega^{CK}_1$, {\rm e.g.,} $\mbox{{\rm T=KP}}\omega$ , {\rm then} $A^{L_\Ome }=_{df}A $. {\rm Define the}
$\Pi_{2}^{\Ome}$-ordinal $|\mbox{{\rm T}}|_{\Pi_{2}^{\Ome}}$ of  {\rm T by}
\[|\mbox{{\rm T}}|_{\Pi_{2}^{\Ome}}=_{df}\inf\{\alpha\leq\omega^{CK}_1 :\forall\Pi_2\mbox{ {\rm sentence }}A(\mbox{{\rm T}}\vdash A^{L_\Ome} \: \Rightarrow \: L_\alpha\models A)\}<\omega^{CK}_1\]
\end{definition}
Here note that $|\mbox{T}|_{\Pi_{2}^{\Ome}}<\omega^{CK}_1$ since we have for any $\Pi_2$ sentence $A$, $\mbox{T}\vdash A^{L_\Ome}  \Rightarrow  L_\Ome\models A$
and $\Ome=\ome^{CK}_1$ is recursively regular, i.e., $\Pi_2$-reflecting.

G. J\"ager \cite{J} shows that $|\mbox{KP}\omega|_{\Pi_{2}^{\Ome}}=\psi_{\Ome}\veps_{\Ome+1}=d_{\Ome}\veps_{\Ome+1}=$Howard ordinal and G. J\"ager and W. Pohlers \cite{J-P} gives the ordinal $|\mbox{KPi}|_{\Pi_{2}^{\Ome}}=\psi_{\Ome}\veps_{I+1}$, where KPi denotes a set theory for recursively inaccessible universes and $I$ the first (recursively) weakly inaccessible ordinal. These include and imply proof-theoretic ordinals of subsystems of second order arithmetic corresponding to set theories. 

$\mbox{KP}\ome$ includes $\mbox{ID}_{1}$: Using (\ref{eqn:BK}), the axiom schema $n\in W(\prec)(\Lrarw_{df}Wf[\prec| n]) \rarw TI[\prec| n,F]$ for arithmetical $\prec$ (note that this expresses the well-founded part $W(\prec)$ is the least fixed point of the operator determined by the formula $\fal m\prec n(m\in X)$) is derivable from $\Sig\mbox{-rfl}\,  (\Del_{0}\mbox{-Coll})$ and {\em Foundation axiom schema} $\fal x[\fal y\in x \vphi(y)\rarw \vphi(x)]\rarw \fal x\vphi(x)$.

\bremark
\benu
\item {\rm (J\"ager, \cite{Jaeger84b}) $|\mbox{KP}\ome_{0}|_{\Pi^{1}_{1}}=\veps_{0}$. In $\mbox{KP}\ome_{0}$ Foundation is restricted to sets $x\in a$.
\item (Rathjen, \cite{Rathjen92}) For $n\geq 2$ $|\Pi_{n}\mbox{-Fund}|_{\Pi_{2}^{\Ome}}=d_{\Ome}\Ome_{n-1}(\ome)$ with $\Ome_{0}(\ome)=\ome\spand \Ome_{n+1}(\ome)=\Ome^{\Ome_{n}(\ome)}$. In $\Pi_{n}\mbox{-Fund}$ Foundation is restricted to $\Pi_{n}$-formulae $\vphi(x)$.

Note that $d_{\Ome}\Ome^{\ome}$ is the Ackermann ordinal, {\it cf}. \cite{R-W93}, and $d_{\Ome}\Ome^{2}=\Gam_{0}$, the first strongly critical number.
}
\eenu
\eremark
{\bf Ramification, level and hierarchy}. In the proof-theoretic analysis of predicative parts of second order arithmetic (Sch\"utte et.al) second order variables$X$ is {\em stratified} into {\em ramified analytic hierarchy} according to contexts (occurrences of $X$ in proofs): $\ome$-models ${\cal M}_{\alp}=(\ome,M_{\alp};0,+,\cdot,\ldots)$. Put $M_{0}=\Del^{0}_{1}$ and let $M_{\alp+1}$ denote the collection of definable subsets of $\ome$ in ${\cal M}_{\alp}$. E.g., $M_{1}=\Pi^{1}_{0}$. Alternatively we can set $M_{\alp}$ as the jump hierarchy. 

For ID theories by Pohlers (local predicativity) the least fixed point $I=I_{<\Ome}$ is stratified into $\bigcup\{I_{\alp}:\alp<\Ome\}$.

In J\"ager's case $L_{\alp}$ do the same job.

\smallskip

KPi is a constructive ZF in a sense: KPi is equivalent to each one of Feferman's $T_{0}$, Martinl\"of's type theory (1984), $\Del^{1}_{2}\mbox{-CA+BI}$. Using the following lemma we see that $\Del^{1}_{2}\mbox{-CA}$ is derived from $\Del_{1}\mbox{-Sep}$.
$Ad(d)$ designates that $d$ is admissible. Note that there is a $\Pi^{0}_{3}$ sentence $\theta$ so that for any transitive $d$, $Ad(d)\lrarw \theta^{d}$, a $\Del_{0}$ formula.

\blem\label{lem:Quant}Let $\sig$ be a limit of admissible ordinals.
\benu
\item For each $\Pi^{1}_{1}$ formula $A(n,X)$ there exists a $\Sig_{1}$ formula $A_{\Sig}(n,X)$ in the language of set theory so that (cf. (\ref{eqn:BK}).)
\[L_{\sig}\models Ad(d)\spand n\in\ome\spand X\incl\ome\spand X\in d\rarw [A(n,X)\lrarw A_{\Sig}^{d}(n,X)]\]
\item
For each $\Sig^{1}_{2}$ formula $F(n,Y)$ with a set parameter $Y$ there exists a $\Sig_{1}$ formula $A_{\Sig}(n,Y)$ so that for
\begin{eqnarray}
&& F_{\Sig}(n,Y)\Lrarw_{df}\exi d[Ad(d)\spand Y\in d\spand A^{d}_{\Sig}(n,Y)] \nonumber \\
&& L_{\sig}\models n\in\ome\spand Y\incl\ome\rarw \{F(n,Y)\lrarw F_{\Sig}(n,Y)\} \label{eq:Quant}
\end{eqnarray}
\item 
For each $\Sig^{1}_{m+1}$ formula $F(n,Y)$ with a set parameter $Y$ there exists a $\Sig_{m}$ formula $F_{\Sig_{m}}(n,Y)$ so that 
\[L_{\sig}\models n\in\ome\spand Y\incl\ome\rarw \{F(n,Y)\lrarw F_{\Sig_{m}}(n,Y)\} \]
\eenu
\elem

\section{Prehistory to Mahlo}\label{sec:5}
J\"ager \cite{Jaeger84a}, Pohlers \cite{Pohlers87}, Sch\"utte \cite{Schuette88}(1984-1988) investigated 
$\rho$-inaccessible ordinals. $0$-inaccessibles are regular cardinals. $(\rho+1)$-inaccessibles are regular fixed points of the function $\pi_{\rho}(\alp)$. For limit $\lam$ $\lam$-inaccessibles are $\rho$-inaccessibles for any $\rho<\lam$. $\pi_{\rho}(\alp)$ is the enumerating function of the continuous closure of $\rho$-inaccessibles. E.g., $\pi_{0}(\alp)=\ome_{\alp}$, $1$-inaccessibles are weakly inaccessibles and $\pi_{1}(\alp)$ are weakly inaccessibles and their limits.

This hierarchy $\pi_{\rho}(\alp)$ of functions reminds us Veblen function $\vphi_{\alp}\bet$.

\section{Recursive notation systems of ordinals}\label{sec:6}
{\em Ordinal diagrams} by Takeuti and us are just finite sequences of symbols together with order relation between them. There may be given set-theoretic interpretations for construtors of o.d's a posteriori. The order relation and constructors on o.d.'s reflect rewriting steps on finite proof figures. To show the well-foundedness of o.d's is the central matter.

While recursive notation systems of ordinals by Buchholz, Rathjen et.al are built in set-theory. First (large) cardinals are supposed to exist, {\it cf}. the subsubsection \ref{subsubsec:8.1.1}. Then define some functions ({\em collapsing functions}) on ordinals to get a structure $(T;<,\Ome,\psi_{\Ome},\ldots)$. Thus we have set-theoretic interpretation and the well-foundedness of the structure in hand a priori assuming the existence of relevant large cardinals. After that the structure is shown to be isomorphic to a {\em recursive} structure $(T;<,\Ome,\psi_{\Ome},\ldots)\simeq (\hat{T};\hat{<},\hat{\Ome},\hat{\psi_{\Ome}},\ldots)$. Further if the latter is shown to be well-founded in a relevant theory, then the assumption of the existence of large cardinals is finally discarded as a figure of speech. 

Another route to dicarding the assumption is to show that either the {\em recursive analogue} of large cardinal suffices to model the structure, \\
$(T;<,\Ome,\psi_{\Ome},\ldots)\simeq (\check{T};<,\check{\Ome}=\ome^{CK}_{1},\check{\psi_{\Ome}},\ldots)$ ({\it cf}. Pohlers \cite{LNM1407} (1989).), or the construction of the structure $(T;<,\Ome,\psi_{\Ome},\ldots)$ is carried (mimiced) in a constructive set theory or a type theory (Rathjen, Griffor, Setzer). When the latter route is pursued, we have to show further that, e.g., a constructive set theory is reduced to a recursive analogue. 

\section{Proof theory of recursively large ordinals}\label{sec:7}
Let ${\cal L}_0$ denote the first order language whose constants are;
$=$(equal), $<$(less than), $0$(zero), $1$(one), $+$(plus), $\cdot$(times), $j $(G\"odel's pairing function on $Ord$),\\
$( )_0 ,( )_1$(projections, i.e., inverses to $j$).

For each $\Delta_0$ formula  ${\cal A}(X,a,b)$ with a binary predicate $X$ in $\calL_{0}\cup\{X\}$ we introduce a binary predicate constant $R^{\cal A}$ and a ternary one $R^{\cal A}_{<}$ by a transfinite recursion on ordinals $a$:
\[
b\in R^{\cal A}_{a} \: \Leftrightarrow_{df} \:  R^{\cal A} (a,b)\:  \Leftrightarrow  \: {\cal A}(R^{\cal A}_{<a} ,a,b) 
\]
  with $R^{\cal A}_{<a} =\sum_{x<a}R^{\cal A}_x = \{(x,y):x<a \, \& \, y\in R^{\cal A}_x\}$.

The language $\calL_1$ is obtained from $\calL_0$ by adding the predicate constants $R^{\cal A}$ and $R^{\cal A}_{<}$ for each bounded formula ${\cal A}(X,a,b)$ in $\calL_{0}\cup\{X\}$.

Let $F:Ord\rightarrow L$ denote (a variant of) the G\"odel's onto map from the class $Ord$ of ordinals to the class $L$ of constructible sets.

The language $\calL_{1}$ is chosen so that the set-theoretic membership relation $\in$ on $L$ is interpretable by a $\Del_{0}$-formula $\in (E,a,b)$ in $\calL_{1}$:
 
\[a\varepsilon b  \Leftrightarrow_{df}F(\alp)\in F(\bet)\Lrarw  \in (E,a,b) \spand  a\equiv b \Leftrightarrow_{df} F(\alp)=F(\bet)\Lrarw =(E,a,b)\]

Thus instead of developing an ordinal analysis of a set theory we can equally develop a proof theory for theories of ordinals.

Every multiplicative principal number $\alp=\ome^{\ome^{\bet}}$ is closed under each function constant in $\calL_{0}$. In particular $\alp$ is closed under the pairing function $j$ and hence each finite sequence $\bar{\bet}<\alp$ is coded by a single $\bet<\alp$. Let $\alp=\langle \alp;0,1,+,\cdot,\ldots,R^{\cal A}|\alp,\ldots\rangle$ denote the $\calL_{1}$-model with the universe $\alp$. We sometimes identify the set $L_{\alp}$ with a multiplicative principal number $\alp$ since $L_{\alp}=F"\alp$.

$\Pi_2^\Ome$-ordinal $|\mbox{T}|_{\Pi_{2}^{\Ome}}$ of a sound and recursive theory T of ordinals is defined similarly as before.

In order to get an upper bound for the $\Pi_2^\Ome$-ordinal $|\mbox{T}|_{\Pi_{2}^{\Ome}}$ of a theory T we attach a {\it term\/} $o(\Gam;P)$ to each sequent $\Gam$ occurring in a proof $P$ in the theory T, which ends with a $\Pi_{2}^{\Ome}$ sentence. The term $o(\Gam;P)$ is built up from atomic diagrams and {\it variables\/} by applying constructors in a system $(O(\mbox{T}),<)$ of o.d.'s for T. Variables occurring in the term $o(\Gam;P)$ are eigenvariables occurring below $\Gam$. Thus the term $o(\Gam_{end};P)$ attched to the endsequent of $P$ is a closed term, i.e., denotes an o.d. Also each redex in our transformation is on the main branch, i.e., the rightmost branch of a proof tree and is the lowermost one. Therefore when we resolve an inference rule $J$ no free variable occurs below $J$.

 Finally set
\[o(P)=d_{\Ome}o(\Gam_{end};P)\in O(\mbox{T})\!|\!\Omega(=\{\alpha\in O(T):\alpha<\Omega\}),\]
where $d_{\Ome}\alp$ is a collapsing function
\[d_{\Ome}:\alp\mapsto d_{\Ome}\alp<\Ome\]

Applied constructors in building the term $o(\Gam;P)$ correspond the inference rules occurring above $\Gam$. For example at an inference rule $(b\exi)$
\[\infer[(b\exi)]{\Gam,\exi x<tA(x)}{\Gam,s<t & \Gam,A(s)}\]
we set with a complexity measure $gr(A)$ of formulae $A$
\[o(\Gam,\exi x<tA(x))=o(\Gam,s<t)\# o(\Gam,A(s))\#s\# gr(A(s))\]
 Note that the instance term $s$ may contain variables, e.g., $s\equiv y\cdot z$. Also at an inference rule $(b\fal)$
\[\infer[(b\fal)]{\Gam,\fal x<t A(x)}{\Gam,x\not<t,A(x)}\]
we substitute the term $t$ for the eigenvariable $x$ in the term $o(\Gam,\fal x<t A(x))$; 
\[o(\Gam,\fal x<t A(x))=o(\Gam,x\not<t,A(x))[x:=t]\]
Also, for example, to analyze (the inference rule corresponding to) the following axiom saying $\Ome$ is $\Pi_{2}$-reflecting
\[\fal u<\Ome [A^{\Ome}(u)\rarw  \exi z<\Ome (u<z\spand A^{z}(u)]\: (A \mbox{ is a } \Pi_{2}\mbox{ formula)}\]
we introduce a new rule together with a new constructor $(\Ome,\alp)\mapsto d_{\Ome}\alp<\Ome$ of o.d.'s:

\[\infer[(c)^{\Ome}_{d_{\Ome}\alp}]{\Gam^{\Ome},A^{d_{\Ome}\alp}}{\Gam^{\Ome},A^{\Ome}}\]
with a set $\Gam$ of $\Sig_{1}$ sentences. $\alp$ is chosen so that $\alp=o(\Gam^{\Ome},A^{\Ome})$.

Now our theorem for an upper bound is stated as follows.
\bth\label{th:up}
If $P$ is a proof of a $\Pi_{2}^{\Ome}$-sentence $A^{\Ome}$ in {\rm T}, then $A^\alpha$ is true with $\alp=o(P)$.
\eth

\section{Reflecting ordinals}\label{sec:8}

\begin{definition} {\rm (Richter and Aczel \cite{R-A}) }
{\rm Let} $X\subseteq Ord$ {\rm denote a class of ordinals and} $\Phi$ {\rm a set of formulae in the language of set theory (or the language of theories of ordinals). Put} $X\!|\!\alpha=_{df}\{\beta\in X :\beta<\alpha \}$. {\rm We say that an ordinal} $\alpha\in Ord$ {\rm is} $\Phi$-reflecting on $X$  {\rm if} 
\[\forall A \in\Phi \mbox{ {\rm with parameters from }} L_\alpha [L_\alpha\models A \: \Rightarrow \: 
\exists\beta\in X\!\!\alpha (L_\beta\models A)]\]
 {\rm If a parameter} $\gamma<\alpha$ {\rm occurs in} $A ${\rm , then it should be understood that} $\gamma<\beta${\rm .}
\\
$\alp$ {\rm is} $\Phi$-reflecting {\rm if} $\alpha$ {\rm is} $\Phi${\rm -reflecting on the class of ordinals} $Ord$.
\end{definition}
This is known to be a recursive analogue to indescribable cardinal $\kap$:
\[\fal R\incl V_{\kap}[\langle V_{\kap},\in,R\rangle\models A\Rarw \exi\alp<\kap(\langle V_{\alp},\in,R\cap V_{\alp}\rangle\models A)]\]
{\bf Facts and definitions}. \cite{R-A}
\begin{enumerate}
\item $\alpha\in Ad \, \& \, \alpha>\omega \: \Leftrightarrow \: \alpha$ is recursively regular $\Leftrightarrow \: \alpha$ is $\Pi_2$-reflecting (on $Ord$) \\
with $Ad=_{df}$ the class of admissible ordinals
\item $\alpha$ is recursively Mahlo $\Leftrightarrow$ $\alpha$ is $\Pi_2$-reflecting on $Ad.$
\item Put $M_n (X)=_{df}\{\alpha\in X:\alpha$ is $\Pi_n$-reflecting on $X\}$. Then for $n>0$,
\[M_{n+1}(Ad)\subseteq M_n^\triangle (Ad), (M_n^\triangle)^\triangle (Ad), etc.,\] 
where $M_n^\triangle$ denotes the diagonal intersection of the operation 
\\
$X \mapsto M_n(X)$.\\
The least $\Pi_{n+1}$-reflecting ordinal is greater than, e.g., the least ordinal in $M_n^\triangle (Ad)$.
\end{enumerate}

From \cite{R-A} we know that $\Pi_{3}$-reflecting ordinals are recursive analogues to $\Pi^{1}_{1}$-indescribable cardinals, i.e., weakly compact cardinals. We say that $\kap$ is $2${\it -regular} if for every $\kap$-bounded $F:{}^{\kap}\kap\rarw{}^{\kap}\kap$ there exists an $\alp$ such that $0<\alp<\kap$ and for any $f\in{}^{\kap}\kap$, if $\alp$ is closed under $f$, then $\alp$ is also closed under $F(f)$. Here $F$ is $\kap$-bounded if
\[\fal f\in{}^{\kap}\kap\fal\xi<\kap\exi\gam<\kap\fal g\in{}^{\kap}\kap[g\gam=f\gam\rarw F(f)(\xi)=F(g)(\xi)]\]
Then $\kap$ is $2$-regular iff $\kap$ is weakly compact.

Let $\kap$ be an admissible ordinal and $\xi<\kap$. We say $\{\xi\}_{\kap}$ maps $\kap$-recursive functions to $\kap$-recursive functions if
\[\fal\bet<\kap[\{\bet\}_{\kap}:\kap\rarw\kap\Rarw \{\{\xi\}_{\kap}(\bet)\}_{\kap}:\kap\rarw\kap]\]
An admissible $\kap$ is said to be $2${\it -admissible} iff for any $\xi<\kap$ if $\{\xi\}_{\kap}$ maps $\kap$-recursive functions to $\kap$-recursive functions, then there exists an $\eta$ such that $\xi<\eta<\kap$ and $\{\xi\}_{\eta}$ maps $\eta$-recursive functions to $\eta$-recursive functions. Then $\kap$ is $2$-admissible iff $\kap$ is $\Pi_{3}$-reflecting. 

\subsection{$\Pi_{2}$-reflection}\label{subsec:8.1}

\subsubsection{A system $O(\Ome)$ of ordinal diagrams}\label{subsubsec:8.1.1}
We define a system $O(\Ome)$ of ordinal diagrams. $O(\Ome)$ is equivalent to Takeuti's system $O(2,1)$ and the Howard ordinal is denoted by the o.d. $d_{\Ome}\veps_{\Ome+1}$.

 Let $0,\Ome,+,\ome^{\alp}$(exponential with base $\ome$) and $d$ be distinct symbols. Each element called ordinal diagram in the set $O(\Ome)$ is a finite sequence of these symbols.
 
$0,\Ome$ are atomic diagrams and constructors in the system $O(\Ome)$ are $+,\ome^{\alp}$ and $d_{\Ome}:\alp\mapsto d_{\Ome}\alp$.\footnote{$\alp$ in $d_{\Ome}\alp$ is not restricted to the case $\alp\geq\Ome$.} Each diagram of the form $d_{\Ome}\alp$ and $\Ome$ are defined to be epsilon numbers:
\[\bet<d_{\Ome}\alp\Rarw\ome^{\bet}<d_{\Ome}\alp\]
The order relations between epsilon numbers are defined as follows.
\benu
\item $d_{\Ome}\alp<\Ome$
\item $d_{\Ome}\alp<d_{\Ome}\bet$ holds if one of the following conditions is fulfilled.
 \benu
 \item $d_{\Ome}\alp\leq K_{\Ome}\bet(\Lrarw_{df}\exi\del\in K_{\Ome}\bet(d_{\Ome}\alp\leq\del))$
 \item $K_{\Ome}\alp<d_{\Ome}\bet(\Lrarw_{df}\fal\gam\in K_{\Ome}\alp(\gam<d_{\Ome}\bet))\spand \alp<\bet$
 \eenu
\item $K_{\Ome}\alp$ denotes the finite set of subdiagrams of $\alp$ which are in the form $d_{\Ome}\gam$, i.e., $K_{\Ome}\alp$ consists of the epsilon numbers below $\Ome$ which are needed for the unique representation of $\alp$ in Cantor normal form.
\eenu 
Then we have the following facts.
\bdes
\item [($<1$)] $d_{\Ome}\alp<\Ome$
\item [($<2$)]$K_{\Ome}\alp<d_{\Ome}\alp$
\item [($<3$)]$K_{\Ome}\alp\leq\alp$
\item[($<4$)] $\bet<\Ome \, \&\, K_{\Ome}\bet<d_{\Ome}\alp \: \Rarw \: \bet<d_{\Ome}\alp$
\edes

An essentially or a collapsibly less than relation $\alp\ll\bet$ is defined by
\[\alp\ll\bet\Lrarw K_{\Ome}\alp<d_{\Ome}\bet\spand\alp<\bet\Lrarw d_{\Ome}\alp<d_{\Ome}\bet\spand \alp<\bet\]

The sytem $O(\Ome)$ is nothing but the notation system $D(\veps_{\Ome+1})$ defined in \cite{R-W93}. Put
\[k_{\Ome}\alp=\max(K_{\Ome}\alp\cup\{0\}) \spand \Ome=\ome_{1}\mbox{(the first uncountable cardinal)}\]
Define sets $D(\alp)$ and ordinals $d_{\Ome}\alp$ by simultaneous recursion on $\alp$ as follows:
\benu
\item $\{\Ome\}\cup (k_{\Ome}\alp+1)\incl D(\alp)$
\item $D(\alp)$ is closed under $+,\ome^{\bet}$.
\item $\del\in D(\alp)\cap\alp\Rarw d_{\Ome}\del\in D(\alp)$
\item $d_{\Ome}\alp=\min\{\xi:\xi\not\in D(\alp)\}$
\eenu
Then we see
\benu
\item $d_{\Ome}\alp<\Ome=\ome_{1}$
\item $d_{\Ome}\bet\leq K_{\Ome}\alp\Rarw d_{\Ome}\bet<d_{\Ome}\alp$
\item $\alp<\bet\spand K_{\Ome}\alp<d_{\Ome}\bet\Rarw d_{\Ome}\alp<d_{\Ome}\bet$
\item $d_{\Ome}\alp=d_{\Ome}\bet\Rarw\alp=\bet$
\item $d_{\Ome}\alp=D(\alp)\cap\Ome$
\item $\alp\in D(\bet)\Lrarw K_{\Ome}\alp<d_{\Ome}\bet$
\eenu

\subsubsection{Finitary analysis}\label{subsubsec:8.1.2}
We explain our approach to an ordinal analysis by taking theories of $\Pi_{2}$ reflecting ordinals as an example.

The fact that $\Ome$ is $\Pi_{2}$ reflecting is expressed by the following inference rule:
\[\infer[(\Pi_{2}\mbox{-rfl})]{\Gam}
{ \Gam,A^{\Ome} & \neg\exi z(t<z<\Ome\wedge A^{z}),\Gam}
\]
for any $\Pi_{2}$-formula $A^{\Ome}\equiv A\equiv\fal x\exi y B(x,y,t)$ with a {\it parameter term} $t$.
$\mbox{T}_{2}$ denotes the theory obtained from $\mbox{T}_{0}$ by adding the inference rule 
$(\Pi_{2}\mbox{-rfl})$.   $\mbox{T}_2$ is formulated in Tait's logic calculus. 

Let $\calL_{c}$ denote the extended language of $\calL_{1}$ obtained by adding an individual constant $\bet$ for each o.d. $\bet<\Ome$.
\[\calL_{c}=\calL_{1}\cup\{\bet\in O(\Ome):\bet<\Ome\}\]
We show
\bth\label{th:T2}
\[\forall\Pi_2\: A(\mbox{{\rm T}}_2\vdash A^{\Ome} \: \Rightarrow \: \exists\alpha\in O(\Omega)\!|\!d_{\Ome}\veps_{\Ome+1}\:A^\alpha ).\]
\eth

 Let $P$ be a proof ending with a $\Pi^{\Ome}_{2}$ sentence $A^{\Ome}$. To each sequent $\Gam$ in $P$, we assign a term $o(\Gam;P)\in{\cal F}$ so that $A ^\alp$ is true with $\alp=d_{\Ome}\alp_0$ and $\alp_0=o(P)$. This is proved by induction on $\alp$.

To deal with the rule $(\Pi_2\mbox{-rfl})$ we introduce a new rule:

\[\infer[(c)^{\Ome}_{d_{\Ome}\alp}]{\Gam,A^{d_{\Ome}\alp}}{\Gam,A^{\Ome}}\]
where $\Gam\subset\Sig^{\Ome}_{1}$ sentences, $A^{\Ome}\equiv\fal x\exi yB$ is a $\Pi^{\Ome}_{2}$-sentence and the following condition have to be enjoyed:
\beqn\label{eq:cnd}
o(\Gam,A^{\Ome})\ll\alp
\eeqn

This rule is plausible in view of the Collapsing Lemma \ref{lem:Collapsing}. 

\blem\label{lem:Collapsing}{\rm (\cite{J})}
Collapsing Lemma: $\vdash^\alpha_\Omega\Gamma \, \& \, \Gamma\subset\Sigma_1 \: \Rightarrow \:  d_{\Ome}\alpha\models\Gamma$
\elem
 where $\beta\models\Gamma \: \Leftrightarrow_{df} \: \bigvee\Gamma^\beta=\bigvee\{\exists x_1<\beta B_1,\ldots,\exists x_n<\beta B_n\}$ ($B_1,\ldots,B_n$ are bounded) is true in the model 
$\langle O(\Omega)\!|\!\beta;+,\cdot,j,\ldots,R^{\cal A}|\beta,\ldots\rangle$.

When a $(\Pi_2\mbox{-rfl})$ is to be analyzed,
\[\infer[(\Pi_{2}\mbox{-rfl})]{\Gam}
{ \Gam,A^{\Ome} & \neg\exi z(t<z<\Ome\wedge A^{z}),\Gam}
\]
roughly speaking, we set $\alp=o(\Gam,A^{\Ome})$ and substitute $d_{\Ome}\alp$ for the variable $z$ [originally $z$ is replaced by $\Ome$], and replace the $(\Pi_2\mbox{-rfl})$ by a $(cut)$.

The inference rule $(\Pi_{2}\mbox{-rfl})$ is resolved as follows:
\[
\infer[J]{\Lam}
{
 \infer[(c)^{\Ome}_{d_{\Ome}\alp}]{\Lam,A^{d_{\Ome}\alp}}
 {
  \infer*{\Lam,A^{\Ome}}{\Gam,A^{\Ome}}
 }
&
 \infer*{\neg A^{d_{\Ome}\alp},\Lam}
 {
  \infer*{}
  {
   \infer*{\neg A^{d_{\Ome}\alp},\Gam}{\del}
  }
 }
}
\]
where 
\benu
\item $\alp=o(\Lam,A^{\Ome})$. 
\item $(c)^{\Ome}_{d_{\Ome}\alp}$ is the new inference rule, which says, if $\Pi^{\Ome}_{2}$-sentence $A^{\Ome}$ is derivable with a $\Sig^{\Ome}_{1}$ side formulae $\Lam$ and an o.d. $\alp$, then we have $\Lam,A^{d_{\Ome}\alp}$, viz. after substituting any $\del<d_{\Ome}\alp$ coming from the right upper part of the $(cut)\, J$ for the universal quantifier $\fal x<\Ome$ in $A^{\Ome}$, we should have $\bet<d_{\Ome}\alp$ for any instance term $\bet<\Ome$ of the existential quantifier $\exi y<\Ome$ in $A^{\Ome}$. 
\item The right upper part of $J$ is obtained by inversion, i.e., substituting the individual constant $d_{\Ome}\alp$ for the variable $z$. $t<d_{\Ome}\alp<\Ome$ follows from $t<\Ome$ and the fact that $t$ is contained in $\alp$, {\it cf}. ($<4$).
\eenu

Then the points are that we have to retain the condition (\ref{eq:cnd}) $o(\Gam,A)\ll\alp$ in the rule $(c)$ and if we have 

\[\infer[(c)]{\Gam,\exi y<d_{\Ome}\alp B}{\infer[(\exi)]{\Gam,\exi yB(y)}{\Gam,B(\bet_{0})}}\]

then it should be the case $\bet_{0}<d_{\Ome}\alp$, i.e., $d_{\Ome}\bet\in K_{\Ome}\bet_{0}\Rarw d_{\Ome}\bet<d_{\Ome}\alp$.
\smallskip
First of all, $d_{\Ome}\bet$ occurs in a proof only because $d_{\Ome}\bet$ was generated at a $(c)$ and then substituted at a $(\Pi_2\mbox{-rfl})$. The latter condition $d_{\Ome}\bet<d_{\Ome}\alp$ is ensured by the former (\ref{eq:cnd}) since $d_{\Ome}\bet\ll o(\Gam,A)\ll\alp$. The former condition (\ref{eq:cnd}) is retained since the only unbounded universal quantifier in $\Gam,A$ is the outermost one $\fal x$ in $A$ and the o.d.$\geq d_{\Ome}\alp$ is forbidden to be substituted for $x$ by the restriction $\exi x<d_{\Ome}\alp$ in $\neg A$.

Observe that there exists a {\em gap} $[d_{\Ome}\alp,\Ome)$ for o.d.'s occurring above a rule $(c)^{\Ome}_{d_{\Ome}\alp}$. Namely if $\bet<\Ome$ occurs above $(c)^{\Ome}_{d_{\Ome}\alp}$, then $\bet<d_{\Ome}\alp$. This follows from the condition (\ref{eq:cnd}) and the fact ($<4$): 
\[\bet<\Ome \, \&\, K_{\Ome}\bet<d_{\Ome}\alp \: \Rarw \: \bet<d_{\Ome}\alp\]

Thus the Theorem \ref{th:T2} was shown by a finitary analysis.

\subsection{Summary of results}\label{subsec:8.2}

\begin{tabular}{|c|c|c|c|}
\hline
ordinal & set-ordinal  & arithmetic & ordinal diagrams \\ 
& theory & & \\
\hline
rec. regular & $KP\ome$ & $\exi{\cal O}, \, ID_1$ & $O(\Ome)^*$   \\
\hline
rec. inacc. & $KPi$  & $\Sig^1_2-AC+BI$,  &  $O(1;I)^*$  \\
& & $SBL$ & \\
\hline
rec. Mahlo & $KPM$ &  & $O(\mu)^*$  \\
\hline
$\Pi_n$-reflecting  & $T_{n}$ & & $O(\Pi_n)$ \\
\hline
$\Pi^1_1$-reflecting & $T_{1}^{1},S(2;1,1)$  &  & $O(2;1,1)$ \\
\hline
\end{tabular}

$^*$ designates that the o.d.'s are shown to be optimal. 

In a letter \cite{W2} A. Weiermann informed me that an inspection of his work in \cite{W1} yields an embedding of $O(\mu)$ in the notation system $T(M)$ by Rathjen \cite{Rathjen90}. Thus via Rathjen's well-ordering proof in \cite{Rathjen94a} we get indirectly that $O(\mu)$ is best possible.

Recently we showed that $\mbox{KPM}\vdash Wo[O(\mu)|\alp]$ for each $\alp<d_{\Ome}\veps_{\mu+1}$ without referring \cite{Rathjen94a}.

\section{Stability}\label{sec:9}
Rreflecting ordinals are too small to model the axiom $\Sig^{1}_{2}\mbox{-CA}$ of second order arithmetic and hence theories for these ordinals are intermediate stages towards $\Sig^{1}_{2}\mbox{-CA}$. We have to consider theories for ordinals below which there are stable ordinals.

\bdf\label{df:10nstbl}
{\rm Let} $\kap$ {\rm and} $\sig<\kap$ {\rm be ordinals and} $k$ {\rm a positive integer. We say that} $\sig$ {\rm is} $(\kap,k)$-stable {\rm if}
\[L_{\sig}\prec_{\Sig_{k}}L_{\kap},\]
{\rm that is, for any} $\Sig_{k}$ {\rm formula} $A$ {\rm with parameters from} $L_{\sig}$ 
\[L_{\kap}\models A \Rarw L_{\sig}\models A.\]
\edf
Note that $(\kap,1)$-stability is equivalent to $\kap$-stability.
\smlskp
{\bf Facts}. (cf.\cite{R-A} and \cite{M}.) For a countable $\sig$,
\benu
\item $\sig$ is $\Pi^1_{0}$-reflecting $\Lrarw \: \sig$ is weakly stable, $\bet$-stable for some $\bet>\alp$.
\item $\sig$ is $\Pi^1_1$-reflecting $\Lrarw \: \sig$ is $\sig^+$-stable.
\item $\Pi^1_1$ on $L_\sig=$inductive on $L_\sig=\Sig_1$ on $L_{\sig^+}$. \\
($\sig^+$ denotes the next admissible to $\sig$.)
\eenu

\subsection{Summary of results}\label{subsec:9.1}

The reason for this turning to stability is that $\Sig^1_{k+1}$-Comprehension Axiom for $k\geq 1$ is interpretable in a universe $L_{\kap}$ such that $L_{\kap}$ has $(\kap,k)$-stable ordinals and $L_{\kap}$ is a limit of admissible sets..

Let $\sig_{0}$ denote a $\Pi_{3}$ sentence in the language of set theory so that a transitive set $x$ is admissible iff $x\models\sig_{0}$, {\it cf}. pp.315-316 in \cite{R-A}. Let $st_{0}(x)$ denote the $\Pi_{0}$ formula:
\[st_{0}(x)\equiv\sig_{0}^{x}\spand x \mbox{ is transitive.}\]
Also for $k\geq 1$ let $st_{k}(x)$ denote a $\Pi_{k}$ formula such that for any admissible $\kap$
\[L_{\kap}\models st_{k}(\sig) \Lrarw L_{\sig}\prec_{\Sig_{k}}L_{\kap}\]
Let $\mbox{KP}\ell^{r}_{k}$ denote a set theory for limits of ordinals $\sig$ with $st_{k}(\sig)$.\footnote{The superscript $r$ in $\mbox{KP}\ell^{r}_{k}$ indicates that the foundation schema is restricted to sets.}
\[(Lim)_{k} \: \fal x\exi y(x\in y\spand st_{k}(y))\]

Using Lemma \ref{lem:Quant} one can model the axiom $\Sig^{1-}_{k+1}\mbox{-CA}$: \\
$\exi X[X=\{n\in\ome: F(n)\}]$ ($F(n)$ is a $\Sig^{1-}_{k+1}$ formula without set parameter.)
in the universe $L_{\kap}$ which contains a $(\kap,k)$-stable ordinal $\sig<\kap$ and $L_{\kap}\models (Lim)_{0}$:
\beqnarrs
\{n\in\ome:L_{\kap}\models F^{set}(n)\} & = & \{n\in\ome:L_{\kap}\models F_{\Sig_{k}}(n)\} \\
& = & \{n\in\ome:L_{\sig}\models F_{\Sig_{k}}(n)\}\in L_{\sig+1}\incl L_{\kap}
\eeqnarrs

For $\Sig_{k+1}$ formula $\vphi(x)\equiv\exi y\theta(y,x)$ and $L_{\kap}\models (Lim)_{k}$
\[\vphi(x)\lrarw \exi\alp[st_{k}(\alp)\spand x\in L_{\alp}\spand\vphi^{L_{\alp}}(x)]\]
This enables us to iterate $\Sig_{1}$-stability proof theory in analysing $\Sig_{k+1}$-stability.

\begin{tabular}{|c|c|c|c|}
\hline
$A+1$ stables & $S(2;A+1)$ & $\Sig^{1-}_2\mbox{-CA}_{1+A+1}$ & $O(2;A+1)^*$ \\
\hline
limit $A$ stables & $S(2;A)$ & $\Sig^1_2\mbox{-AC +BI}+$ & $O(2;A)^*$  \\
 & & $\Sig^{1-}_2\mbox{-CA}_{A}$ & \\
\hline
$<\ome$-stables  &  $S(2;<\ome)$ & $\Sig^1_2\mbox{-CA}_0$  & $O(2;<\ome)^*$  \\
\hline
 $\ome$-stables, & $\Sig_1$-Sep, & $\Sig^1_2\mbox{-CA+BI}$ & $O(2;\ome)^*$  \\
  nonprojectible & $S(2;\ome)$ & & \\
\hline
 $<\ome^\ome$-stables & $S(2;<\ome^{\ome})$ & $\Sig^1_3\mbox{-DC}_0$  & $O(2;<\ome^\ome)^*$ \\ 
\hline
 $<\veps_0$-stables &  $S(2;<\veps_{0})$ & $\Sig^1_3\mbox{-DC}$ &  $O(2;<\veps_0)^*$ \\
\hline
 $\Pi_2(St)$-reflecting  & $\Pi_1$-Coll., &  $\Sig^1_3\mbox{-AC+BI}$ & $O(2;I)^*$ \\
 on stables $St$ &  $S(2;I)$ & & \\
\hline 
$A+1$ 2-stables & $S(3;A+1)$ & $\Sig^{1-}_{3}\mbox{-CA}_{1+A+1}$ & $O(3;A+1)^*$ \\
&&& (?) \\
\hline
$<\ome$ 2-stables  &  $S(3;<\ome)$ & $\Sig^1_{3}\mbox{-CA}_0$  & $O(3;<\ome)^*$  \\
&&& (?) \\
\hline
 $<\ome^\ome$ 2-stables &  $S(3;<\ome^{\ome})$ & $\Sig^1_{4}\mbox{-DC}_0$  &  $O(3;<\ome^{\ome})^*$ \\ 
 &&& (?) \\
\hline
 $<\veps_0$ 2-stables &   $S(3;<\veps_0)$ & $\Sig^1_{4}\mbox{-DC}$ & $O(3;<\veps_0)^*$  \\
 &&& (?) \\
\hline
\end{tabular}
\smlskp
\[\Sig^{1-}_2\mbox{-CA}_{1+A+1}\,: \: \exi\{X_a\}_{a<_{A}1\oplus A}\fal a<_{A}1\oplus A(X_a=\{n: F(n,a,X_{<_{A}\,a})\})\]

 $S(2;I)$ denotes a theory of ordinals $I$ such that $I$ is $\Pi_2(St)$-reflecting, where $St$ denotes the set of stable ordinals below $I$ and $\Pi_2(St)$ the set of $\Pi_2$ formulae $A$ in the language $\calL_{1}\cup\{St\}$ so that the predicate constant $St$ may occur. Then the set theory $\mbox{KP}\ome+\Pi_1 \mbox{-Collection+V=L}$ is interpretable in $S(2;I)$.

\subsection{Proof theory for $\Pi^{1}_{1}$-reflection}\label{subsec:9.2}
A baby case for ordinals below which there is a stable ordinal is an ordinal $\pi^{+}$ such that $\pi^{+}$ is the next admissible to a $\pi^{+}$-stable ordinal $\pi$, viz. $\Pi^{1}_{1}$-reflecting ordinal.
Such a universe $L_{\pi^{+}}$ can be modelled in a theory $T^{1}_{1}$ for positive elementary inductive definitions on $L_{\pi}$:
 Fix an $X$-positive formula $A\equiv A(X^+,a)$ in the language 
$\calL_{1}\cup\{X\}$. Let $Mp$ denote the set of multiplicative principal numbers $a\leq\pi$. Define a ternary predicate $I_<$ by:
\[\fal a\in Mp\fal b<a^+[I^a_{<b}=\bigcup_{d<b}I^a_d=\bigcup_{d<b}\{c<a:A^a(I^a_{<d},c)\}]\]

That is to say, for each $a\in Mp,a\leq\pi$ and $b<a^+$,  $I^a_{<b}$ is the inductively generated subset of $a=\{c:c<a\}$ by the positive formula $A$ on the model $<a;+,\cdot,\ldots,R^{\cal A},\ldots>$, uniformly with respect to the multiplicative principal number $a$.

Then the axioms of the theory $T^{1}_{1}$ say that the universe $\pi^+$ is $\Pi_2$-reflecting and the axiom $(\Pi^1_1\mbox{-rfl})$:
\[\fal c<\pi[c\in I^\pi_{<\pi^+}\rarw \exi \bet\in (c,\pi)\cap Mp(c\in I^\bet_{<\bet^+})].\]
 where $c\in I^a_{<a^+}\Lrarw_{df}\exi z<a^+\, A^a(I^a_{<z},c)$.

The theory $T^{1}_{1}$ is designed so that a theory $S_{1}^{1}$ for ordinals $\pi^{+}$ with $L_{\pi}\prec_{\Sig_{1}}L_{\pi^{+}}$ is interpretable in $T_{1}^{1}$.

Let us examine the crucial case.
\[\infer[(\Pi^1_1\mbox{-rfl})\, J]{}
{
 \neg(\alp<b<\pi),\fal x<b^+\neg A^b(I^b_{<x},\alp) 
& 
 \infer[(\exi)]{\exi x<\pi^+ A^\pi(I^\pi_{<x},\alp)}
              {A^\pi(I^\pi_{<\xi},\alp)}
}\]
with $\alp\in I^\pi_{<\pi^+}\equiv\exi x<\pi^+ A^\pi(I^\pi_{<x},\alp)$, etc.

First consider the easy case:\\
{\bf Case1}. $\xi<\pi$: Then the above $(\Pi^1_1\mbox{-rfl})\, J$ says that $\pi$ is $\Pi_\xi$-reflecting. So define $\sig=d_\pi$ such that $\xi,\alp<\sig<\pi$ and substitute $\sig$ for the variable $b$.\\
Second the general case:\\
{\bf Case2}. $\xi\geq\pi$: Pick a $\sig=d_\pi$ as above and substitute $\sig$ for $b$. We need to compute a $\xi'$ such that $\sig\leq\xi'<\sig^+$ and resolve the $(\Pi^1_1\mbox{-rfl})\, J$:

\[\infer[(cut)]{}
{
 \neg A^\sig(I^\sig_{<\xi'},\alp)
&
 \infer[(c)^\pi_\sig \, I]{A^\sig (I^\sig _{<\xi'},\alp)}
  {
    \infer[J]{A^\pi(I^\pi_{<\xi},\alp)}
            {
              \alp\not\in I^b_{<b^+} 
            &
              \alp\in I^\pi_{<\pi^+},A^\pi(I^\pi_{<\xi},\alp)
            }
   }
}
\]

The problem is that we have to be consistent with the part

\[\infer{}
{
 \neg A^\sig(I^\sig_{<\xi'},\alp)
&
 \infer{A^\sig (I^\sig _{<\xi'},\alp)}
  {A^\pi(I^\pi_{<\xi},\alp)}
}
\]  
namely
\[A^\sig(I^\sig_{<\xi'},\alp)\lrarw A^\pi(I^\pi_{<\xi},\alp)\]

This requires a function $F :\xi\mapsto\xi'$ such that
\bdes
\item[$(F1)$] $F$ is order preserving, \\
and in view of {\bf Case1},
\item[$(F2)$] $F$ is identity on$<\pi$, i.e., $\xi\in dom(F)|\pi\Rarw F(\xi)=\xi$
\item [$(F3)$]$rng(F)<\sig^+$.
\edes
Note that, here, $dom(F)$ is a {\em proper subset} of $\{\xi\in O(\Pi^1_1):\xi<\pi^+\}$ with a system $O(\Pi^1_1)$ of o.d.'s for the theory $T^1_1$. We can safely set
    \[dom(F )=\{\xi\in O(\Pi^1_1):K_\pi\xi<\sig\} \]
 i.e., subdiagram $\bet<\pi$ in $\xi\in dom(F)$ is $<\sig$ since $dom(F)$ is the set of o.d.'s that may occur in the upperpart of the $(c)^\pi_\sig\, I$. Especially we have
    \[dom(F )|\pi=O(\Pi^1_1)|\sig\]
This would be possible since there exists a {\em gap} $[\sig,\pi)$ for o.d.'s occurring above the rule $(c)^\pi_\sig$.

   Can we take the function $F$ as a collapsing function, e.g., $d_\pi$? The answer is no. We cannot expect for $\xi,\zeta\in dom(F)$, that $\xi<\zeta \Rarw \xi\ll_\pi\zeta$ or something like an essentially less than relation. And what is worse is that the function $F$ have to preserve atomic sentences in $\calL_0$.
\bdes
\item[$(F4)$] $F$ preserves atomic sentences in $\calL_0$, i.e., diagrams of $\calL_0$ models 
 \\
     $< dom(F ) ; +,\cdot,\ldots>$ and $< rng(F ) ; +,\cdot,\ldots>$. 
\edes

To sum up $(F1)-(F4)$,
\bdes
\item[(*)] $F$ is an embedding from $\calL_0$ models $<dom(F);+,\cdot,\ldots>$ \\
to $<rng(F); +,\cdot,\ldots>$ over $O(\Pi_1^1)|\sig$.
\edes

Now our solution for $F$ is a trite one: a {\em substitution} $[\pi:=\sig]$.
\bdes
\item [$(F5)$]$F(\xi)=\xi$ if $\xi<\pi\, (\Lrarw \xi<\sig)$
\item [$(F6)$]$F$ commutes with $+$ and the Veblen function $\vphi$, e.g., \\
$F(\xi+\zeta)=F(\xi)+F(\zeta)$.
\item [$(F7)$]$F(\pi)=\sig$ and $F(\pi^+)=\sig^+$.
\item [$(F8)$] $F(d_{\pi^+}\bet)=d_{\sig^+}F(\bet)$.
\edes
Assume $\pi<\xi<\pi^+$ with a strongly critical $\xi$. Such a $\xi$ is of the form $d_{\pi^+}\bet$ and is introduced when a $(\Pi_2\mbox{-rfl})$ for the universe $\pi^+$ is resolved. Then this $F$ meets (*), i.e., $(F1)$: Note that we have
\bdes
\item[$(F9)$] $F(K_{\pi^+}\bet)=K_{\sig^+}F(\bet)$,
\edes
and by definition $d_{\pi^+}\bet<d_{\pi^+}\gam \, \Lrarw \, 1. \, \bet<\gam \,\&\, K_{\pi^+}\bet< d_{\pi^+}\gam \mbox{ or } 2. \, d_{\pi^+}\bet\leq K_{\pi^+}\gam$ 
and similarly for $\sig^+$.

In fact a miniature $[\sig,\veps_{\sig^{+}+1})$ of $[\pi,\veps_{\pi^{+}+1})$ is formed by a realisation $F$ of the Mostowski collapsing function. 
\smallskip

In this way we can resolve a $(\Pi_1^1\mbox{-rfl})$ by setting $\xi'=F(\xi)$: each o.d. $\xi$ in the uppersequent of a $(c)$ is replaced by $F(\xi)$ in the lowersequent.

\section{The future:uncountable cardinals}\label{sec:10}
It is important to find an {\it equivalent and right axiom} in begining proof-theoretic analysis for recursively large ordinals. For example nonprojectible ordinal $\kap$ was analysed by us as a limit of $\kap$-stable ordinals. At least for us the latter formulation was essential: if we adopted other axioms, e.g., there is no $\kap$-recursive injection $f:\kap\rarw\alp<\kap$ or $L_{\kap}\models \Sig_{1}\mbox{-Separation}$, then an analysis of these axioms would be difficult for us. Therefore in this final section we give an equivalent condition for $\kap$ to be an uncountable cardinal. The condition remains a submodel condition saying $\kap$ has an appropriate submodel. So it may be possible to analyse such a universe by extending proof theory for $\Sig_{k}$-stability in the near future.

\bdf\label{df:crd} {\rm Let} $\sig$ {\rm be a recursively regular ordinal and} $\ome\leq\alp<\kap<\sig$.
\benu
\item {\rm We say that} $\kap<\sig$ {\rm is} a $\sig$-cardinal{\rm , denoted} $L_{\sig}\models \kap \mbox{ is a cardinal }$ iff 
\[L_{\sig}\models\fal\alp\in[\ome,\kap)[\mbox{there is no surjective map } f:\alp\rarw\kap]\]
\item 
\[L_{\sig}\models card(\alp)<card(\kap)\Lrarw_{df}L_{\sig}\models \mbox{there is no surjective map } f:\alp\rarw\kap\]
\eenu
\edf

\bth\label{th:crdlocal}Let $\sig$ be a recursively regular ordinal and $\kap,\alp$  multiplicative principal numbers with $\ome\leq\alp<\kap<\sig$. Then the following conditions are mutually equivalent:
\benu
\item 
\beqnarr
&& \exi(\pi,\pi\kap,\pi\sig)[\alp<\pi\leq\pi\kap<\pi\sig<\kap\spand  \nonumber \\
&& \fal \Sig_{1}\,\vphi\fal a<\pi(L_{\sig}\models\vphi[\kap,a]\rarw L_{\pi\sig}\models\vphi[\pi\kap,a])] \label{eqn:crdlocal1}
\eeqnarr
\item 
\beqn\label{eqn:crdlocal2}
{\cal P}(\alp)\cap L_{\sig}\incl L_{\kap}
\eeqn
\item 
\beqn\label{eqn:crdlocal3}
L_{\sig}\models card(\alp)<card(\kap)
\eeqn
\eenu
\eth

In what follows $\sig$ denotes a recursively regular ordinal and $\alp,\kap$ multiplicative principal numbers with $\ome\leq\alp<\kap<\sig$.

\blem\label{lem:crdlocal1}
(\ref{eqn:crdlocal1})$\Rarw$(\ref{eqn:crdlocal2})
\elem
\bprf First note that $L_{\kap}\prec_{\Sig_{1}}L_{\sig}$. Define a $\Del_{1}$-partial map $S:dom(S)\rarw{\cal P}(\alp)\cap L_{\kap}\, (dom(S)\incl\kap)$ as follows. First set $S_{0}=\emptyset$ and let $S_{\bet}$ denote the $<_{L}$ least $X\in{\cal P}(\alp)\cap L_{\kap}$ such that $\fal\gam<\bet(X\not\in S_{\gam})$. 

It suffices to show that ${\cal P}(\alp)\cap L_{\sig}\incl\{S_{\bet}\}=rng(S)$. Suppose there exists an $X\in{\cal P}(\alp)\cap L_{\sig}$ so that $\fal\bet<\kap(S_{\bet}\neq X)$ and let $X_{0}$ denote the $<_{L}$-least such set. Then $X_{0}$ is $\Sig_{1}$ definable in $L_{\sig}$: there exists a $\Del_{1}$ formula
\[\vphi(X,\alp,\kap)\Lrarw_{df} \theta(X,\alp,\kap)\spand \fal Y<_{L}X\neg\theta(Y,\alp,\kap)\]
with
\[\theta(X,\alp,\kap)\Lrarw_{df}X\incl\alp\spand\fal\bet<\kap(S_{\bet}\neq X)\]
so that
\[L_{\sig}\models\vphi(X_{0},\alp,\kap)\spand L_{\sig}\models\exi !X\vphi(X,\alp,\kap)\]
By (\ref{eqn:crdlocal1}) we have $L_{\pi\sig}\models\exi X\vphi(X,\alp,\pi\kap)$, i.e., there exists the $<_{L}$-least $X_{1}\in{\cal P}(\alp)\cap L_{\pi\sig}\incl {\cal P}(\alp)\cap L_{\kap}$ such that $\fal\bet<\pi\kap(<\kap)(S_{\bet}\neq X_{1})$. This means that $X_{1}=S_{\pi\kap}$. We show $X_{1}=X_{0}$. This yields a contradiction. 

Denote $x\in a$ by $x\in^{+}a$ and $x\not\in a$ by $x\in^{-}a$. For any $\gam<\alp$, again by (\ref{eqn:crdlocal1}) we have
\beqnarrs
&& \gam\in^{\pm} X_{0}\Lrarw L_{\sig}\models\exi X(\gam\in^{\pm}X\spand \vphi(X,\alp,\kap)) \Rarw \\
&& L_{\pi\sig}\models\exi X(\gam\in^{\pm}X\spand \vphi(X,\alp,\pi\kap)) \Lrarw \gam\in^{\pm}X_{1}
\eeqnarrs
\eprf

\blem\label{lem:crdlocal2}
(\ref{eqn:crdlocal2})$\Rarw$(\ref{eqn:crdlocal3})
\elem
\bprf Argue in $L_{\sig}$. Suppose there exists a surjective map $f:\alp\rarw\kap$. Pick a surjective map (in $L_{\sig}$) $g:\kap\rarw L_{\kap}$. Let $F:\alp\rarw{\cal P}(\alp)(\cap L_{\sig})$ denote the map given by
\[F(\bet)=\left\{
\begin{array}{ll}
g(f(\bet)) & g(f(\bet))\in{\cal P}(\alp) \\
\emptyset & \mbox{otherwise}
\end{array}
\right.
\]
Then by (\ref{eqn:crdlocal2}), ${\cal P}(\alp)\cap L_{\sig}\incl L_{\kap}$ $F$ is surjective. Also $F\incl\alp\times L_{\kap}$ is $\Del_{0}$ and hence $F\in L_{\sig}$ by $\Del_{0}$-Separation. Define $X\in{\cal P}(\alp)\cap L_{\sig}$ by
\[X=\{\bet<\alp:\bet\not\in F(\bet)\}\]
Then $X=F(\gam)$ for some $\gam<\alp$ and $\gam\in X\Lrarw \gam\not\in F(\gam)=X$. This is a contradiction.
\eprf

\blem\label{lem:crdlocal3}
(\ref{eqn:crdlocal3})$\Rarw$(\ref{eqn:crdlocal1})
\elem
\bprf Since $\alp$ is a multiplicative principal number, each finite sequence $\bar{\bet}<\alp$ is coded by a single $\bet<\alp$.

We define a $\Sig_{1}$ subset $X$ of $L_{\sig}$ ($\Sig_{1}$-Skolem hull of $\alp\cup\{\alp,\kap\}$ in $L_{\sig}$): Let $\{\vphi_{i}:i\in\ome\}$ denote an enumeration of $\Sig_{1}$-formulae of the form $\vphi_{i}\equiv\exi y\theta_{i}(x,y;z,u,v)$ with a fixed variables $x,y,z,u,v$. Set for $\bet<\alp$
\[r(i,\bet)\simeq \mbox{ the } <_{L} \mbox{ least } c\in L_{\sig}[L_{\sig}\models\theta_{i}((c)_{0},(c)_{1};\bet,\alp,\kap)]\]
\[h(i,\bet)\simeq(r(i,\bet))_{0}\]
and
\[X=rng(h)=\{h(i,\bet)\in L_{\sig}:i\in\ome, \bet<\alp\}\]
Clearly $r$ and $h$ are partial $\Sig_{1}$ map whose domains are $\Sig_{1}$ subset of $\ome\times\alp$.
First note that
\beqn\label{eqn:crdlocal3.1}
\alp\cup\{\alp,\kap\}\incl X
\eeqn
Next we show
\bclm\label{clm:crdlocal3.1} For any $\Sig_{1}(X)$-sentence $\vphi(\bar{a})$ with parameters $\bar{a}$ from $X$
\[L_{\sig}\models \vphi(\bar{a})\Lrarw X\models \vphi(\bar{a})\]
Namely
\[X\prec_{\Sig_{1}}L_{\sig}\]
\eclm
{\bf Proof} of Claim \ref{clm:crdlocal3.1}. Suppose $L_{\sig}\models\exi v\theta(v,\bar{a})$ with $\bar{a}\incl X$. It suffices to show that there exists a $b\in X$ so that $L_{\sig}\models \theta(b,\bar{a})$. For each $a_{k}\in\bar{a}$ pick a $\Sig_{1}$-formula $\vphi_{i_{k}}\equiv \exi y\theta_{i_{k}}(x,y;z,u,v)$ and $\bet_{k}<\alp$ so that $h(i_{k},\bet_{k})\simeq a_{k}$. Then 
\[L_{\sig}\models \exi v\exi\bar{x}[\theta(v,(\bar{x})_{0})\spand \bigwedge[x_{k} \mbox{ is the } <_{L} \mbox{ least } w \theta_{i_{k}}((w)_{0},(w)_{1};\bet_{k},\alp,\kap)]\]
where $(\bar{x})_{0}=(x_{0})_{0},\ldots,(x_{n})_{0}$ with $\bar{x}=x_{0},\ldots,x_{n}$.
 Hence the assertion follows.
\eprfclm

Suppose for the moment that the $\Sig_{1}$-subset $dom(h)\incl\ome\times\alp$ is an element of $L_{\sig}$ ($\sig$-finite). Then $h$ is $\Del_{1}$ and $X=rng(h)$ is $\Del_{1}$-subset of $L_{\sig}$. We show
\bclm\label{clm:crdlocal3.2} Assume $dom(h)\in L_{\sig}$. Then there exist a triple $(\pi,\pi\kap,\pi\sig)$ satisfying (\ref{eqn:crdlocal1}).
\eclm
{\bf Proof} of Claim \ref{clm:crdlocal3.2}. By Claim \ref{clm:crdlocal3.1} and the Condensation Lemma ({\it cf}. p.80 in \cite{Devlin}.) we have an isomorphism (Mostowski collapsing function) $F:X\lrarw L_{\pi\sig}$ for an ordinal $\pi\sig\leq\sig$ such that $F| Y=id| Y$ for any transitive $Y\incl X$.

We show first that $\pi\sig<\kap$. Suppose $\kap\leq\pi\sig$.
The collapsing function $F$ is defined by the following recursion:
\[F(x)=\{F(y):y\in x\spand y\in X\}\]
Since $X$ is $\Del_{1}$, $F$ is a $\Del_{1}$-function. The $\Del_{1}$ map $h$ maps $\ome\times\alp$ onto $X$ and the $\Del_{1}$ map $F$ maps $X$ onto $L_{\pi\sig}\supseteq L_{\kap}$. Hence the composition $F\circ h$ maps $\ome\times\alp$ maps onto $L_{\pi\sig}$. Let $G$ denote a restriction of $F\circ h$ so that $rng(G)=L_{\kap}$. Then its domain $dom(G)$ is a $\Del_{1}$-subset of $dom(h)$ and hence $dom(G)\in L_{\sig}$. Therefore by combining a surjective map from $\alp$ onto $\ome\times\alp$ we get a $\Del_{1}$ map $f\incl\alp\times\kap$ such that $dom(f)=\alp$ and $rng(f)=\kap$. $\Del_{1}$-Separation in $L_{\sig}$ yields $f\in L_{\sig}$. This is a contradiction since $card(\alp)<card(\kap)$ in $L_{\sig}$. Thus we have shown $\pi\sig<\kap$.

 Let $\pi$ denote the least ordinal not in $X$ and set $\pi\kap=F(\kap)$. Then $F(a)=a$ for any $a<\pi$. Also clearly $\alp<\pi\leq\pi\kap<\pi\sig<\kap$. For a $\Sig_{1}$ sentence $\vphi[\kap,a]$ with a parameter $a<\pi$ assume $L_{\sig}\models\vphi[\kap,a]$. Then $X\models\vphi[\kap,a]$ and hence $L_{\pi\sig}\models\vphi[\pi\kap,a]$ as desired. 
\eprfclm
\smlskp
Thus it remains to show the 
\bclm\label{clm:crdlocal3.3} $dom(h)\in L_{\sig}$.
\eclm
{\bf Proof} of Claim \ref{clm:crdlocal3.3}. $dom(h)=\{(i,\bet)\in\ome\times\alp: L_{\sig}\models\exi c\theta_{i}((c)_{0},(c)_{1};\bet,\alp,\kap)\}$. Let $\sig^{*}$ denote the $\Sig_{1}$-projectum of $\sig$. $dom(h)$ is a $\Sig_{1}$-subset of $\ome\times\alp\lrarw\alp$. Thus it suffices to show ({\it cf}. Theorem 6.11\footnote{Any $\sig$-r.e. subset of $\bet<\sig^{*}$ is $\sig$-finite for admissible $\sig$.} on p.177, \cite{Barwise}.)
\[\alp<\sig^{*}\]
Suppose $\sig^{*}\leq\alp$. Let $F:\sig\rarw\sig^{*}$ denote a $\Sig_{1}$ injection and $f=F|\kap$ the restriction of $F$ to $\kap$. Then $f\in L_{\sig}$ would be an injection from $\kap$ to $\sig^{*}\leq\alp$. This is a contradiction since 
\[L_{\sig}\models card(\alp)<card(\kap)\lrarw \mbox{there is no injective map } f:\kap\rarw\alp\]
\eprf

\bth\label{th:crd}Let $\sig$ be a recursively regular ordinal and $\kap$ a multiplicative principal number with $\ome<\kap<\sig$. Then the following conditions are mutually equivalent:
\benu
\item 
\beqnarrs
&& \fal\alp\in[\ome,\kap)\exi(\pi,\pi\kap,\pi\sig)[\alp<\pi\leq\pi\kap<\pi\sig<\kap\spand  \\
&& \fal \Sig_{1}\,\vphi\fal a<\pi(L_{\sig}\models\vphi[\kap,a]\rarw L_{\pi\sig}\models\vphi[\pi\kap,a])] 
\eeqnarrs
\item 
\[
\fal\alp\in[\ome,\kap)[{\cal P}(\alp)\cap L_{\sig}\incl L_{\kap}]
\]
\item 
\[
L_{\sig}\models \kap \mbox{ is a cardinal}>\ome
\]
\eenu
\eth

\bth\label{th:crdexi}Let $\sig$ be a recursively regular ordinal and $\gam$ a multiplicative principal number with $\gam<\sig$. The following conditions are mutually equivalent.
\benu
\item $\exi\kap\in Mp|\sig \exi(\pi,\pi\kap,\pi\sig)[\gam<\pi\leq\pi\kap<\pi\sig<\kap\spand 
 \fal \Sig_{1}\,\vphi\fal a<\pi(L_{\sig}\models\vphi[\kap,a]\rarw L_{\pi\sig}\models\vphi[\pi\kap,a])] $
\item ${\cal P}(\gam)\cap L_{\sig}\in L_{\sig}$
\item $L_{\sig}\models\exi\kap(\kap \mbox{ is a cardinal}>\gam)$
\eenu
\eth
\bprf By Theorem \ref{th:crdlocal} it suffices to show the last condition assuming the second one. Assume ${\cal P}(\gam)\cap L_{\sig}\in L_{\sig}$. Then there exists a multiplicative principal $\kap_{0}<\sig$ such that ${\cal P}(\gam)\cap L_{\sig}\incl L_{\kap_{0}}$. By Lemma \ref{lem:crdlocal2} we have $L_{\sig}\models card(\gam)<card(\kap_{0})$. Let $\kap<\sig$ denote the least ordinal satisfying this. Then we claim that $L_{\sig}\models \kap \mbox{ is a cardinal}$. For suppose $L_{\sig}\not\models card(\alp)<card(\kap)$ for some $\alp$ with $\gam<\alp<\kap$ by $L_{\sig}\models card(\gam)<card(\kap)$. Pick a surjective map $f\in L_{\sig}$ with $f:\alp\rarw\kap$. Also pick a surjective map $g\in L_{\sig}$ with $g:\gam\rarw\alp$ by the minimality of $\kap$, i.e., $L_{\sig}\not\models card(\gam)<card(\alp)$. The composition $f\circ g:\gam\rarw \kap$ is a surjective map in $L_{\sig}$ contrdicting $L_{\sig}\models card(\gam)<card(\kap)$.
\eprf

\end{document}